\documentstyle[12pt]{article}

\newtheorem{theorem}{Theorem}
\newtheorem{proposition}{Proposition}
\newcommand{\QED}{{\hfill$\Box$\medskip}}

\newtheorem{corollary}{Corollary}

\begin{document}
\begin{center}
{\Large \bf Computations of the Yamabe invariant}\\
\vspace{1 cm}
Jimmy Petean\\
Max-Planck Institut f\"ur Mathematik,\\
Bonn, Germany\\
\end{center}

\vspace{1 cm}
\begin{abstract}
For a compact connected manifold $M$ of dimension  $n\geq 4$ 
with no metric of positive
scalar curvature, we prove that the Yamabe invariant  is unchanged
under surgery on spheres of dimension different from 1, $n-2$ and $n-1$. 
We use this result to 
give new exact computations  of the Yamabe invariant in dimension four and
display new examples of compact four-manifolds
which do not admit Einstein metrics. 
\end{abstract}

\section{Introduction}

Consider a smooth compact connected manifold $M^n$ of dimension $n$.
It is well known that set of the critical points of the functional 

$${\cal S} (g) =
\frac{\int\limits_M s_g \  dvol_g}{(Vol_g (M))^{\frac{n-2}{n}}}$$

\noindent 
(on the space ${\cal M}$ of all Riemannian metrics on $M$)
is the set of  Einstein metrics on $M$. But it is also well known that this
functional is never bounded below or above and hence it is not possible
to prove existence of critical points by minimizing or maximizing the
functional.
An alternative strategy to find critical points is a minimax argument,   
the first ideas of which were introduced by 
Yamabe  in \cite{Yamabe}. First we have to restrict
the functional to some fixed conformal class; the functional is bounded
below in this subspace and we can define the Yamabe constant of the
conformal class  ${\cal C}$ by

$$Y(M, {\cal C})= \inf\limits_{g\in {\cal C}} {\cal S}  (g).$$

Moreover, the  infimum is always achieved by a smooth metric, which
therefore has constant scalar curvature.
Yamabe claimed to prove this result in \cite{Yamabe}, but his  proof
contained an error. This was recognized by Trudinger (see \cite{Trudinger}),
and the problem was  finally solved in several steps by Trudinger 
\cite{Trudinger},
Aubin \cite{Aubin} and  Schoen \cite{Schoen2}. 

\vspace{.3cm}

We now define the {\it Yamabe invariant} of $M$ by

$$Y(M)=\sup\limits_{\cal{C} }Y(M,\cal{C})$$

\noindent
(where the supremum is taken over the family of all conformal classes of
metrics on $M$).

This invariant was introduced by O. Kobayashi
in \cite{Kobayashi} and  it is also 
frequently called the {\it sigma constant} 
of $M$ \cite{Schoen3}.
Note  that the Yamabe invariant of $M$
is an invariant of the smooth structure of $M$.

\vspace{.3cm}

We will be concerned in this paper with manifolds for which the Yamabe 
invariant is non-positive. 
This is equivalent to say that the manifold does not admit any metric
of strictly positive scalar curvature. In this case,  there is always a unique
metric of constant scalar curvature and unit volume  
in each conformal class of metrics on $M$,  and the Yamabe invariant  is the
supremum of the scalar curvature over the space of metrics of constant scalar
curvature and unit volume. Furthermore, 
if the Yamabe invariant is realized by such a metric, then
the metric is a critical point of the functional ${\cal S}$ and  it 
is therefore an Einstein metric. 

This is the minimax approach to find Einstein metrics
we mentioned in the beginning of this introduction.
For the proofs of these last statements, and additional
references, see \cite[Section 1]{Schoen3}.

\vspace{.3cm}

Recall also that, still assuming that $M$ does not admit any metric of
positive scalar curvature,

$$Y(M)=-\left( \inf\limits_{\cal M} {\int}_M |s_g|^{n/2} \ dvol_g 
\right)^{2/n},$$

\noindent
where ${\cal M}$ is, as before, the space of all Riemannian metrics on
$M$ (see \cite{Anderson2, LeBrun}). Hence all our results 
regarding the Yamabe invariant of a manifold can also be 
reinterpreted in terms of the infimum 
over ${\cal M}$ of the
$L^{n/2}$ norm of the scalar curvature.

\vspace{.3cm}

We will prove,

\begin{theorem}: Suppose  $M$ is a connected  smooth compact manifold
of dimension $n\geq  4$ with $Y(M)\leq 0$. 
Let $\widehat{M}$ be a manifold obtained
from $M$ by performing surgeries on spheres of dimension different from
1, n-1 and n-2. Then $Y(\widehat{M} )=Y(M)$.
\end{theorem}

The proof of Theorem 1 will be based on the results of \cite{Petean},
where it is proved that surgery in codimension greater than 2 does not
decrease the Yamabe invariant (assuming that the invariant is
non-positive).

\vspace{.3cm}

Some computations of the invariant have been carried out in low dimensions.
In dimension four, LeBrun computed the Yamabe invariant of all compact complex
surfaces of K\"{a}hler type which do not admit metrics of
positive scalar curvature \cite{LeBrun}. 
We  will use Theorem 1 and the computations in \cite{LeBrun}  to give new
exact computations of the Yamabe invariant for 
four-dimensional manifolds. 
Note that the computation of the invariant in dimension two follows from 
the Gauss-Bonnet formula and in dimension three, in the non-positive case,
would follow from Anderson's
program for the hyperbolization conjecture \cite{Anderson}.

We will also show that it follows from Theorem 1 that
the  Yamabe invariant of a manifold of dimension greater than four is 
unchanged  under
connecting sum with the product of two spheres (always assuming that the
invariant is non-positive).

\vspace{.3cm}

Finally, we will use these results to display new examples of 
compact four-dimensional manifolds which 
do not admit Einstein metrics. For instance, we
will prove:

\begin{theorem}: The connected sum of a compact complex hyperbolic 
four-manifold
with any number of copies of $S^1 \times S^3$  does not admit any Einstein
metric.
\end{theorem}

Note that similar results,  in the case of {\it real}  hyperbolic 
four-dimensional manifolds,
have been proved by Sambusetti in \cite{Sambusetti}.

Numerous other examples will be given in the last section of this paper.

\section{Computations of the invariant}

We will now prove  Theorem 1. We will use a ``double  surgery'' argument
and the following result, proved in \cite[Theorem 1]{Petean}:

\begin{proposition}: 
If $\widehat{M}$ is obtained from $M$ by performing surgery
in codimension at least 3 and $Y(M)\leq 0$, then $Y(\widehat{M})\geq Y(M)$.
\end{proposition}

Since we know that, under certain  restrictions in codimension,
surgery can only  increase the Yamabe invariant, what we will do
is to try to kill the surgery we realized by another surgery
that also verifies the restriction in codimension.

\vspace{.2cm}

There is a canonical way to undo a surgery done on a manifold. 
Let ${\cal S}^k$ be a $k$-dimensional sphere embedded in the 
manifold $M$ with trivialized normal bundle.
Let $\widehat{M}$
be the manifold obtained by doing surgery on ${\cal S}^k$
(of course,
$\widehat{M}$ will depend on the homotopy class of the trivialization of the
normal bundle).
The boundary of a tubular neighbourhood $U$ of ${\cal S}^k$ is diffeomorphic
to $S^k \times S^{n-k-1}$  and $\widehat{M}$ is obtained by gluing
$D^{k+1} \times S^{n-k-1}$ to $M-U$ along their boundaries. Doing
 surgery
on the (n-k-1)-{\it belt sphere}, $\{ 0\} \times S^{n-k-1} \subset
D^{k+1} \times S^{n-k-1}$, in $\widehat{M}$, obviously produces a manifold
diffeomorphic to $M$. 

\vspace{.2cm}

Assume now that $k\leq n-3$ and $n-k-1\leq  n-3$. Then it follows  from
Theorem  2 and the previous observations
that $Y(M)\leq Y(\widehat{M})\leq Y(M)$. Therefore
$Y(M)=Y(\widehat{M})$. Note that $\widehat{M}$ does not admit any
metric of positive scalar curvature  because if that were the case
then $M$ would also admit such a metric 
(since the family of manifolds which admit metrics of positive
scalar curvature is closed under surgery on codimension greater
than two from the results of Gromov-Lawson
and Schoen-Yau; see \cite{Gromov, Schoen}).

\vspace{.2cm}

The previous  condition on $k$ can also be written as $2\leq  k\leq n-3$.
To conclude the proof of Theorem 1 we therefore need only to deal with the case
$k=0$.
There is another  canonical way to undo surgery in the 0-dimensional case.
Doing 0-dimensional surgery on a connected manifold $M$ 
means that we take out 2 disjoint balls in $M$
and we glue the boundaries to the boundary of 
$S^{n-1}\times [0,1]$. Pick any point 
$q\in S^{n-1}$ and consider the loop obtained by joinning the endpoints of
$\{ q\} \times [0,1]$ with a smoothly embedded 
curve inside of $M$. The result will be a smoothly embedded circle in
the manifold obtained by doing the 0-dimensional surgery on $M$.
Doing surgery on such
a circle will  undo the 0-dimensional surgery. Another way to see this is
the following: the manifold obtained by doing 0-dimensional surgery on $M$ is 
$M\#(S^1 \times S^{n-1})$, and doing surgery on the circle of 
$S^1 \times S^{n-1}$ we get $S^n$. Hence the result of the two surgeries is
$M\#S^n =M$.

When $n\geq 4$, 
this proves that the Yamabe invariant is also unchanged under 0-dimensional
surgeries and, therefore, we have finished the proof of  Theorem 1.

\QED

\vspace{.4cm}

{\it Remark}: It is
also clear  that the Yamabe invariant is unchanged under 
surgery on an (n-1)-sphere  
$ S$ unless  $S$ separates $M$. 

\vspace{.3cm}

Let $\widehat{M}= M\# (S^k \times S^{n-k})$, where $M$ is a smooth connected
compact manifold of dimension $n\geq 5$. Doing surgery  on any of the spheres
$S^k$ or $S^{n-k}$ in $\widehat{M}$ produces a manifold diffeomorphic to $M$.
And, of course, either $k$ or $n-k$ is different  from 1 and n-2
(note also that $M$ admits a metric of positive scalar
curvature if and only if $\widehat{M}$ does). Therefore,
it follows from Theorem 1 that,

\begin{corollary}: Let $M$ be  any  
compact smooth connected manifold $M$  of dimension $n\geq 5$. Assume that
the Yamabe invariant of $M$ is non-positive. 
Then $Y(M)=Y\left( M\#\ (S^k \times  S^{n-k})\right) $ for all $0\leq k\leq n$.
\end{corollary}

The previous result  is definitely not true in the 4-dimensional
case. Recall that C.T.C. Wall proved in \cite{Wall} that if two compact 
simply connected four
manifolds $M^1$ and $M^2$ have isomorphic intersection
forms, then for $l$ big enough the
connected sum of $M^1$ with $l$ copies of $S^2 \times S^2$ is diffeomorphic
to the connected sum of $M^2$ with $l$ copies of $S^2 \times S^2$. 
Suppose that we have a pair  $M^1$, $M^2$ of 4-manifolds
with isomorphic intersection forms such that 
$Y(M^1 )<0$ while $Y(M^2 )>0$. 
Then it  follows from the previous comments (and the fact that the family of
manifolds admitting positive scalar curvature metrics is closed under 
connected sums \cite{Gromov, Schoen}), that for $l$ big enough,
$Y(M^1 \# l(S^2 \times S^2 ))>0$.  
For example, we can take $M^1$ to be  a smooth hypersurface of 
odd degree $d\geq 5$ in ${\bf CP}^3$. C. LeBrun  proved
in \cite{LeBrun2} that $Y(M^1 )<0$, while it is well known, from
Freedman's classification of compact simply connected topological
four-manifolds,  that $M^1$ is
homeomorphic to a connected sum of copies of ${\bf CP}^2$ (some of these
copies taken with the reversed orientation); and such a connected sum
has  positive Yamabe  invariant (since the Fubini-Study metric in
${\bf CP}^2$ has positive scalar curvature).

\vspace{.3cm}

Note that  given a  
4-manifold $M$ and any  ``trivial $S^1$'' in $M$ (i.e. a 
small, simple loop in a neighbourhood of a point), doing surgery on such
a circle produces $M\# (S^2 \times S^2)$. Hence we can see from the
previous observation that Theorem 1  is
close to optimal:

\begin{proposition}: There are examples of 4-manifolds $M$, with negative
Yamabe invariant, on which one can strictly increase the Yamabe invariant
by doing surgery on dimension 1. It follows that there are also examples
where one can strictly decrease the invariant by doing surgery in dimension 2.
\end{proposition}

We can now give some new exact computations of the 
invariant in dimension 4. Using Seiberg-Witten techniques, C. LeBrun 
proved in \cite{LeBrun2} that if $X$ is a minimal compact complex 
surface of general  type and $M$ is  obtained by
blowing up $X$ any number of times, then
 
$$Y(M)=Y(X)=-4\pi\sqrt{2c_1^2 (X)}.$$

Note also that doing 0-dimensional surgery on a connected four-manifold
$M$ produces a manifold diffeomorphic to $M\# (S^1 \times S^3 )$. Therefore,
it follows from  Theorem 1 that,

\begin{proposition}: Let $M$ be a compact four-dimensional 
manifold with non-positive
Yamabe invariant. Then $Y(M)=Y(M\# (S^1 \times S^3 ))$.
\end{proposition}

Using the results of LeBrun mentioned previously we get new  exact 
computations of the invariant:

\begin{proposition}: Let $X$ be a minimal compact complex surface
of general type. Let $M$ be a manifold obtained by taking the connected sum
of any number of copies of $S^1 \times S^3$ with
the blow-up of $X$ at any number of points. Then
$Y(M)=-4\pi\sqrt{2c_1^2 (X)}$.
\end{proposition}

\section{Manifolds without Einstein metrics}

Finally, the results from the previous section produce 
numerous new examples of compact
four-manifolds
which do not admit any Einstein metric. Let $M$ be obtained by blowing up
$l$ times the minimal compact complex surface of general type $X$. 
Let $N=N(k,l)$ be the manifold obtained as the 
connected sum of $M$ with $k$ copies of $S^1 \times S^3$. Suppose that
$N$ admits an Einstein metric $g$. Then it 
follows from the Gauss-Bonnet formula
that

$$2\chi (N)\pm 3\tau (N)=\frac{1}{4{\pi}^2}
\int\limits_N ( 2|W_{\pm}|^2 +(s_g^2 /24)) \  
dvol_g.$$

But,
$$\frac{1}{4{\pi}^2}
\int\limits_N ( 2|W_{\pm}|^2 +(s_g^2 /24))\   
dvol_g \ \ \geq \ \ 
\frac{1}{4{\pi}^2}
\int\limits_N
\frac{s_g^2}{24}\ 
dvol_g \ \ 
\geq$$

$$\geq \ \ \frac{1}{96{\pi}^2}(Y(N))^2 \ \ =\ \  \frac{1}{3}c_1^2 (X).$$

\vspace{.2cm}

\noindent
where the last inequality follows because the Yamabe invariant of $N$ is
the negative of  
the square root of the infimum over the
space of all Riemannian metrics on $N$ of the 
$L^2$ norm of the scalar curvature (we have already 
mentioned this result in the introduction)
and the final equality is precisely Corollary 2.

\vspace{.3cm}

Hence, we have seen that if $N$  admits an Einstein metric then 

$$2\chi (N)-3|\tau (N)|\geq (1/3) c_1^2 (X).$$ 

This implies in many cases that
such an Einstein metric cannot exist. 

For instance,
let us now prove Theorem 2.  Let $M$ be a (compact)
complex
hyperbolic surface. In this case, one knows that 
$\chi (M)=3\tau (M)>0$; and if $N$ is obtained as the
connected sum of $M$ and $k$ copies of $S^1 \times S^3$, the above 
inequality reads: 
$3\tau (M) -4k \geq 3\tau (M)$. Therefore, $N$ cannot admite any Einstein
metric.

\vspace{.2cm}

{\it Remark}: Note that the Hitchin-Thorpe inequality only assures that
an Einstein metric does not exist on $M\# k(S^1 \times S^3 )$ when $k\geq (3/4)
\tau (M)$.

\vspace{.3cm}

For a more general case, 
let $X$ be a minimal compact complex surface of general type
which has non-positive signature. 
As before, let $M$ be the blow up of $X$ at $l$ points and
let $N$ be the connected sum of $M$ and $k$ copies of $S^1 \times S^3$. 
Note  that
$M$ and $N$ also have 
non-positive signature. If $N$ admits an Einstein metric we
get from the previous computations that $c_1^2 (X)-l -4k\geq (1/3) c_1^2 (X)$.
Therefore,

\begin{theorem}: If $k> (1/6)c_1^2 (X) -l/4$, 
then $N=N(k,l)$ does not
admit an Einstein metric.
\end{theorem}

{\it Remark}: In this case the Hitchin-Thorpe inequality only assures the
non-existence of an Einstein metric when $k\geq (1/4)c_1^2 (X) -l/4$.

\vspace{.3cm}

{\bf Acknowledgements}: The  author would like to thank Claude LeBrun
for his support and for numerous observations related to this work. He  would
also like to thank
Vyacheslav Krushkal for many useful discussions, and  the directors
and staff of the Max-Planck Institut for their hospitality during
the preparation of this work.

\end{document}